\documentclass[12pt,a4paper]{article}
\usepackage{amsfonts}
\usepackage{amssymb}
\usepackage{mathrsfs}
\usepackage{amsmath}
\usepackage{amsthm}
\usepackage{enumerate}
\usepackage{cite}
\usepackage[all]{xy}
\usepackage{extarrows}
\usepackage{indentfirst}
\allowdisplaybreaks
\newtheorem{defn}{Definition}[section]
\newtheorem{thm}[defn]{Theorem}
\newtheorem{lem}[defn]{Lemma}
\newtheorem{prop}[defn]{Proposition}
\newtheorem{cor}[defn]{Corollary}
\newtheorem{eg}[defn]{Example}
\newtheorem{re}[defn]{Remark}

\newcommand\relphantom[1]{\mathrel{\phantom{#1}}}
\newcommand{\bdefn}{\begin{defn}}
\newcommand{\edefn}{\end{defn}}
\newcommand{\bthm}{\begin{thm}}
\newcommand{\ethm}{\end{thm}}
\newcommand{\blem}{\begin{lem}}
\newcommand{\elem}{\end{lem}}
\newcommand{\bprop}{\begin{prop}}
\newcommand{\eprop}{\end{prop}}
\newcommand{\bcor}{\begin{cor}}
\newcommand{\ecor}{\end{cor}}
\newcommand{\beg}{\begin{eg}}
\newcommand{\eeg}{\end{eg}}
\newcommand{\bre}{\begin{re}}
\newcommand{\ere}{\end{re}}
\newcommand{\bpf}{\begin{proof}}
\newcommand{\epf}{\end{proof}}
\newcommand{\benu}{\begin{enumerate}}
\newcommand{\eenu}{\end{enumerate}}
\newcommand{\bc}{\begin{center}}
\newcommand{\ec}{\end{center}}
\newcommand{\bea}{\begin{eqnarray}}
\newcommand{\eea}{\end{eqnarray}}
\newcommand{\Bea}{\begin{eqnarray*}}
\newcommand{\Eea}{\end{eqnarray*}}
\newcommand{\beq}{\begin{equation}}
\newcommand{\eeq}{\end{equation}}
\newcommand{\Beq}{\begin{equation*}}
\newcommand{\Eeq}{\end{equation*}}
\newcommand{\bspl}{\begin{split}}
\newcommand{\espl}{\end{split}}

\begin{document}

\title{\textbf{On the structure of graded Leibniz triple systems}
\author{ Yan Cao$^{1,2},$  Liangyun Chen$^{1}$
 \date{{\small {$^1$  School of Mathematics and Statistics, Northeast Normal
 University,\\
Changchun 130024, China\\
\small {$^2$}Department of Basic
 Education,
 Harbin University of
Science and Technology\\ Rongcheng Campus,  Rongcheng 264300,
China}}}}} \maketitle
\date{}

\begin{abstract}

We study the structure   of a Leibniz triple system $\mathcal{E}$
graded by an arbitrary abelian group $G$ which is considered of
arbitrary dimension and over an arbitrary base field $\mathbb{K}$.
We show that $\mathcal{E}$ is of the form
$\mathcal{E}=U+\sum_{[j]\in \sum^{1}/\sim} I_{[j]}$ with $U$ a
linear subspace of the 1-homogeneous component $\mathcal{E}_{1}$ and
any  ideal  $I_{[j]}$ of $\mathcal{E}$, satisfying
$\{I_{[j]},\mathcal{E},I_{[k]}\}
=\{I_{[j]},I_{[k]},\mathcal{E}\}=\{\mathcal{E},I_{[j]},I_{[k]}\}=0$
if $[j]\neq [k]$, where the
relation $\sim$ in $\sum^{1}=\{g \in G \setminus \{1\} : L_{g}\neq 0\}$, defined by $g \sim h$ if and only if $g$ is connected to $h$.\\

\end{abstract}
\renewcommand{\thefootnote}{\fnsymbol{footnote}}
\footnote[0]{ Corresponding author(L. Chen): chenly640@nenu.edu.cn.}
\footnote[0]{Supported by  NNSF of China (Nos. 11171055 and
11471090) and Scientific
Research Fund of Heilongjiang Provincial Education Department
 (No. 12541184). }

\section{Introduction}
 Leibniz triple systems were introduced by Bremner and S\'{a}nchez-Ortega \cite{BS}. Leibniz triple systems were defined in a functorial manner using
 the Kolesnikov-Pozhidaev algorithm, which took the defining identities for a variety of algebras and produced the defining identities for the corresponding variety
  of dialgebras \cite{K}.  In \cite{BS}, Leibniz triple systems were obtained by applying the Kolesnikov-Pozhidaev algorithm to Lie triple systems. The study of gradings on Lie algebras begins in the 1933 seminal Jordan's work, with the purpose of formalizing Quantum Mechanics\cite{J}. Since then, the interest on gradings on different classes of algebras has been remarkable in the recent years, motivated in part by their application in physics and geometry \cite{BL5238999,BL52567,BL5234,BL52}. Recently, in \cite{BL52,BL528, BL523, L2, AAO, AAO2}, the structure of arbitrary graded Lie algebras, graded Lie superalgebras, graded commutative algebras, graded
  Leibniz algebras and graded  Lie triple systems have been determined by the techniques of connections of roots.
   Our work is essentially motivated by the work on  graded  Lie triple systems\cite{AAO2} and the work on split Leibniz triple systems \cite{AAO2cao}.

Throughout this paper, Leibniz triple systems $\mathcal{E}$ are
considered of arbitrary dimension and over an arbitrary base field
$\mathbb{K}$.  This paper proceeds as follows. In section 2, we
establish the preliminaries on graded Leibniz triple systems theory.
In section 3, we show that such an arbitrary Leibniz triple system
is of the form  $\mathcal{E}=U+\sum_{[j]\in \sum^{1}/\sim} I_{[j]}$
with $U$ a subspace of  $\mathcal{E}_{1}$  and any ideal $I_{[j]}$
of $\mathcal{E}$,  satisfying $\{I_{[j]},\mathcal{E},I_{[k]}\}
=\{I_{[j]},I_{[k]},\mathcal{E}\}=\{\mathcal{E},I_{[j]},I_{[k]}\}=0$
if $[j]\neq [k]$, where the relation $\sim$ in $\sum^{1}$, defined
by $g \sim h$ if and only if $g$ is connected to $h$.

\section{Preliminaries}

\bdefn{\rm\cite{BL5234cao}}  A \textbf{right Leibniz algebra} $L$ is a vector space over a base field $\mathbb{K}$ endowed with a bilinear product
$[\cdot,\cdot]$ satisfying the Leibniz identity
$$[[y, z], x] = [[y, x], z] + [y, [z, x]],$$
for all $x, y, z \in L$.
\edefn

\bdefn{{\rm\cite{BS}}} A \textbf{Leibniz triple system} is a vector space  $\mathcal{E}$ endowed with a trilinear
operation $\{\cdot,\cdot,\cdot\}: \mathcal{E}\times \mathcal{E}\times \mathcal{E}\rightarrow \mathcal{E}$ satisfying
\begin{gather}
\{a,\{b,c,d\},e\}\!=\!\{\{a,b,c\},d,e\}
\!-\!\{\{a,c,b\},d,e\}\!-\!\{\{a,d,b\},c,e\}\!+\!\{\{a,d,c\},b,e\},\label{VIP1}\\
\{a,b,\{c,d,e\}\}\!=\!\{\{a,b,c\},d,e\}\!-\!\{\{a,b,d\},c,e\}\!-\!\{\{a,b,e\},c,d\}\!+\!\{\{a,b,e\},d,c\},\label{VIP2}
\end{gather}
for all $a, b, c, d, e \in \mathcal{E}$.
\edefn

\beg\label{eg}
A Lie triple system gives a Leibniz triple system with the same ternary product.
If $L$ is a Leibniz algebra with product $[\cdot,\cdot]$, then $L$ becomes a Leibniz triple system by putting $\{x,y,z\}=[[x,y],z]$. More examples refer to \cite{BS}.
\eeg

\bdefn{{\rm\cite{BS}}} Let $I$ be a subspace of a Leibniz triple system  $\mathcal{E}$. Then $I$ is called a \textbf{subsystem} of $\mathcal{E}$, if $\{I,I,I\}\subseteq I;$ $I$ is called an \textbf{ideal} of $\mathcal{E}$,
if $\{I,\mathcal{E},\mathcal{E}\}+\{\mathcal{E},I,\mathcal{E}\}+\{\mathcal{E},\mathcal{E},I\}\subseteq I$.
\edefn

\bprop{{\rm\cite{BL5238999cao}}}\label{38888}
Let $\mathcal{E}$ be a  Leibniz triple system. Then the following assertions hold.

$\rm(1)$  $J$ is generated by $\{\{a,b,c\}-\{a,c,b\}+\{b,c,a\}: a,b,c \in \mathcal{E}\}$, then $J$ is an ideal of $\mathcal{E}$  satisfying $\{\mathcal{E},\mathcal{E},J\}=\{\mathcal{E},J,\mathcal{E}\}=0$.

$\rm(2)$  $J$ is generated by $\{\{a,b,c\}-\{a,c,b\}+\{b,c,a\}: a,b,c \in \mathcal{E}\}$, then $\mathcal{E}$ is a Lie triple system if and only if $J=0$.

$\rm(3)$ $\{\{c,d,e\},b,a\}-\{\{c,d,e\},a,b\}-\{\{c,b,a\},d,e\}+\{\{c,a,b\},d,e\}-\{c,\{a,b,d\},e\}-\{c,d,\{a,b,e\}\}=0$, for all $a, b, c, d, e$ $\in$ $\mathcal{E}$.
\eprop

\bdefn{{\rm\cite{BS}}}\label{uni Leib envelop} The \textbf{standard embedding} of a Leibniz triple system $\mathcal{E}$ is the two-graded right Leibniz algebra $L =
L^{0}\oplus L^{1}$, $L^{0}$ being the $\mathbb{K}$-$\rm span$ of $\{x \otimes y,  \ x, y \in \mathcal{E} \}$, $L^{1}: =\mathcal{E}$ and where the product is
given by
$$[(x \otimes y, z), (u \otimes v,w)]:= (\{x, y, u\} \otimes  v - \{x, y, v\} \otimes u
+ z \otimes w, \{x, y, w\} +\{z,u, v\}-\{z,v, u \}).$$
\edefn

Let us observe that $L^{0}$  with the product induced by the one in $L =
L^{0}\oplus L^{1}$  becomes a
 right Leibniz algebra.

\bdefn Let $\mathcal{E}$ be a Leibniz triple system. It is said that $\mathcal{E}$ is graded by means
of an abelian group $G$ if it decomposes as the direct sum of linear subspaces
$$\mathcal{E}=\bigoplus_{g \in G}\mathcal{E}_{g}$$
where the homogeneous
components satisfy $\{\mathcal{E}_{g}, \mathcal{E}_{h}, \mathcal{E}_{k}\} \subset \mathcal{E}_{ghk}$ for any $g, h, k \in G$
$($denoting by juxtaposition the product in $G$$)$. We
call the support of the grading
the set $\sum^{1} := \{g \in G \setminus \{1\}: \mathcal{E}_{g} \neq 0\}.$
\edefn

The usual regularity
conditions will be understood in the graded sense. That is,
a subtriple of $\mathcal{E}$ is a linear subspace $S$ satisfying $\{S, S, S\} \subset S$ and such that splits
as $S = \oplus_{g \in G}S_{g}$
 with any $S_{g} = S \cap \mathcal{E}_{g}$.

Let $L$ be an arbitrary Leibniz algebra over $\mathbb{K}$. As usual, the term grading will always
mean abelian group grading, that is, a decomposition in linear subspaces $L$ =$ \oplus_{g \in G}L_{g}$ where $G$ is  an abelian group and the homogeneous spaces satisfy $[L_{g}, L_{h}] \subset L_{gh}$.
We also
call the support of the grading the set $\{g \in G \setminus \{1\} : L_{g}\neq 0\}$.

\bprop \label{111222}
Let $\mathcal{E}$ be a $G$-graded Leibniz triple system and let $L$ = $L^{0} \oplus L^{1}$ be
its standard embedding algebra, then $L^{0}$ is a $G$-graded Leibniz algebra.
\eprop

\bpf
Define $L^{0}_{1}
:= \sum_{g \in G}[\mathcal{E}_{g}, \mathcal{E}_{g^{-1}}]$  and $L^{0}_{g}
:= \sum_{h \in G}[\mathcal{E}_{h}, \mathcal{E}_{h^{-1}g}]$ for any $g \in G \setminus \{1\}$.
Clearly $L^{0}_{1}+\sum_{g \in G \setminus \{1\}}L^{0}_{g}\subseteq L^{0}$. Conversely, since
 $L^{0}=[\mathcal{E}, \mathcal{E}]=[\bigoplus_{g \in G}\mathcal{E}_{g}, \bigoplus_{h \in G}\mathcal{E}_{h}]\subseteq$
$$L^{0}_{1}+\sum_{g \in G \setminus \{1\}}L^{0}_{g},$$
we get $L^{0}=L^{0}_{1}+\sum_{g \in G \setminus \{1\}}L^{0}_{g}$.

The direct character of the sum can be checked as follows. If $x \in  L^{0}_{g}\cap (\sum_{h \in G \setminus \{g\}} L^{0}_{h})$, then for any $q \in G $ and $y \in \mathcal{E}_{q}$ we have $[x,y] \in \mathcal{E}_{gq}\cap (\sum_{h \in G \setminus \{g\}} \mathcal{E}_{hq})$ and so $[x,y]=0$. From here $[x, \mathcal{E}]=0$ and so $x=0$. Hence we can  write
$$L^{0}=L^{0}_{1}\oplus(\bigoplus_{g \in G \setminus \{1\}}L^{0}_{g}).$$
Finally, we have $$[L^{0}_{g}, L^{0}_{h}]\subseteq L^{0}_{gh}$$
for any $g,h \in G$. Indeed,
$$[L^{0}_{g}, L^{0}_{h}]=\sum_{k,l \in G}[[\mathcal{E}_{k},\mathcal{E}_{k^{-1}g}],[\mathcal{E}_{l},\mathcal{E}_{l^{-1}h}]]\subset$$
$$[\mathcal{E}_{k},[\mathcal{E}_{k^{-1}g},[\mathcal{E}_{l},\mathcal{E}_{l^{-1}h}]]]+[[\mathcal{E}_{k},[ \mathcal{E}_{l},\mathcal{E}_{l^{-1}h}]],\mathcal{E}_{k^{-1}g}]\subset$$
$$[\mathcal{E}_{k},[\mathcal{E}_{k^{-1}g}, L^{0}_{h}]]+[[\mathcal{E}_{k}, L^{0}_{h}], \mathcal{E}_{k^{-1}g}]\subset$$
$$[\mathcal{E}_{k}, \mathcal{E}_{k^{-1}gh}]+[\mathcal{E}_{kh}, \mathcal{E}_{k^{-1}g}]\subset L^{0}_{gh}.$$

Observe that for any $g,h \in G$ we have $[\mathcal{E}_{g},
\mathcal{E}_{h}]\subset L_{gh}^{0}.$ \epf

 In the following, we shall
denote  by  $\sum^{0}$ the support of the graded Leibniz algebra
$L^{0}$.

\section{Connections and gradings}
From now on, $\mathcal{E}$ denotes a graded Leibniz triple system with support $\sum^{1}$, and
$$\mathcal{E} =\bigoplus_{g \in G}\mathcal{E}_{g}=\mathcal{E}_{1}\oplus(\bigoplus_{g \in \sum^{1}}\mathcal{E}_{g})$$
the
corresponding grading. Denote by $-\sum^{i} = \{-g : g \in \sum^{i} \}, i = 0, 1$.

\bdefn Let $g$ and $h$  be two elements in $\sum^{1}$. We say that $g$ is
connected to
$h$ if there exist $g_{1}, g_{2},\cdots,g_{2n+1} \in \pm \sum^{1}\cup \{1\}$ such that

$\rm(1)$ $\{g_{1}, g_{1}g_{2}g_{3},\cdots, g_{1}g_{2}g_{3}\cdots g_{2n}g_{2n+1}\}\subset  \pm \sum^{1}$,

$\rm(2)$  $\{g_{1}g_{2}, g_{1}g_{2}g_{3}g_{4},\cdots, g_{1}g_{2}g_{3}\cdots g_{2n}\}\subset  \pm \sum^{0}$,

$\rm(3)$ $g_{1}=g$ and $g_{1}g_{2}g_{3}\cdots g_{2n}g_{2n+1} \in \{h, h^{-1}\}$.

\noindent We also say that $\{g_{1},\cdots, g_{2n+1}\}$ is a connection from  $g$ to $h$.
\edefn

\bprop \label{6787777}
The relation $\sim$ in $\sum^{1}$, defined by $g \sim h$ if and only if $g$ is connected to $h$ is an equivalence relation.
\eprop
\bpf
The proof is completely analogously to  \cite[Proposition 3.1]{BL523cao}.
\epf

By Proposition \ref{6787777} the connection relation is an equivalence relation in $\sum^{1}$ and
so we
can
consider the quotient set
$\sum^{1}/ \sim= \{[g] : g \in \sum^{1}\}$,
becoming $[g]$ the set of elements in the support of the grading which are
connected
to $g$. By the definition of $\sim$, it is
clear that if $h \in [g]$ and $h^{-1} \in  \sum^{1}$ then $h^{-1} \in [g]$.

Our goal in this section is to associate an adequate subtriple $I_{[g]}$ to any $
[g]$. Fix
$g \in \sum^{1}$, we start by defining
$$\mathcal{E}_{1,[g]} := span_{\mathbb{K}}\{\{\mathcal{E}_{h}, \mathcal{E}_{k}, \mathcal{E}_{(hk)^{-1}}\} : h \in [g], k \in [g] \cup\{1\}\}     \subset \mathcal{E}_{1}$$

\noindent and
$V_{[g]} :=\oplus_{h \in [g]}\mathcal{E}_{h}$. Finally, we denote by $\mathcal{E}_{[g]}$ the direct sum of the two subspaces above, that is,
$$\mathcal{E}_{[g]} := \mathcal{E}_{1,[g]} \oplus V_{[g]}.$$

\bprop \label{caoyan222}
For any $g \in \sum^{1}$, the graded linear subspace $\mathcal{E}_{[g]}$ is a subtriple of
$\mathcal{E}$.
\eprop
\bpf
We have to check that $\mathcal{E}_{[g]}$ satisfies
$$\{\mathcal{E}_{[g]},\mathcal{E}_{[g]},\mathcal{E}_{[g]}\}=\{\mathcal{E}_{1,[g]} \oplus V_{[g]},\mathcal{E}_{1,[g]} \oplus V_{[g]},\mathcal{E}_{1,[g]} \oplus V_{[g]}\}\subset \mathcal{E}_{[g]}.$$

\noindent Since $\mathcal{E}_{1,[g]}\subset \mathcal{E}_{1}$ we clearly have
$$\{\mathcal{E}_{1,[g]},\mathcal{E}_{1,[g]},  V_{[g]}\}+\{\mathcal{E}_{1,[g]}, V_{[g]}, \mathcal{E}_{1,[g]}\}+\{V_{[g]}, \mathcal{E}_{1,[g]},\mathcal{E}_{1,[g]}\}\subset V_{[g]}.$$

\noindent It is easy to verify that $$\{\mathcal{E}_{1,[g]},\mathcal{E}_{1,[g]},\mathcal{E}_{1,[g]}\}\subset \mathcal{E}_{1,[g]}.$$
\noindent Indeed, $\{\mathcal{E}_{1,[g]},\mathcal{E}_{1,[g]},\mathcal{E}_{1,[g]}\}\subset \{\mathcal{E}_{1,[g]},\mathcal{E}_{1},\mathcal{E}_{1}\}\subset \mathcal{E}_{1,[g]}.$
\noindent Moreover, we also have $$\{\mathcal{E}_{1,[g]},V_{[g]},V_{[g]}\}\subset \mathcal{E}_{[g]}.$$
\noindent In fact, if $\{\mathcal{E}_{1}, \mathcal{E}_{h}, \mathcal{E}_{k}\} \neq 0$ for some $h, k \in [g]$ then $h \in \sum^{0}$ and $hk \in \sum^{1}$ $\cup$ $\{1\}$. From here, if $hk \neq 1$ and $\{g_{1},\cdots, g_{2n+1}\}$ is a
connection from $g$ to $h$
then $\{g_{1},\cdots, g_{2n+1}, 1, k\}$ is a
connection from $g$ to $hk$ in
case $g_{1} g_{2}\cdots g_{2n+1} = h$ and  $\{g_{1},\cdots, g_{2n+1}, 1, k^{-1}\}$ in case  $g_{1} g_{2}\cdots g_{2n+1} = h^{-1}$ being so $hk \in [g].$ If $hk=1$ clearly $\{\mathcal{E}_{1}, \mathcal{E}_{h}, \mathcal{E}_{h^{-1}}\}\subset \mathcal{E}_{1,[g]}.$ We have showed $\{\mathcal{E}_{1}, \mathcal{E}_{h}, \mathcal{E}_{k}\}\subset \mathcal{E}_{[g]}$ and so $\{\mathcal{E}_{1,[g]},V_{[g]},V_{[g]}\}\subset \mathcal{E}_{[g]}.$

Next, let us show that $$\{V_{[g]},\mathcal{E}_{1,[g]},V_{[g]}\}\subset \mathcal{E}_{[g]}.$$
\noindent if  $\{\mathcal{E}_{h}, \mathcal{E}_{1},  \mathcal{E}_{k}\} \neq 0$ for some $h, k \in [g]$ then $h \in \sum^{0}$ and $hk \in \sum^{1}$ $\cup$ $\{1\}$. From here, if $hk \neq 1$ and $\{g_{1},\cdots, g_{2n+1}\}$ is a
connection from $g$ to $h$
then $\{g_{1},\cdots, g_{2n+1}, 1, k\}$ is a
connection from $g$ to $hk$ in
case $g_{1} g_{2}\cdots g_{2n+1} = h$ and  $\{g_{1},\cdots, g_{2n+1}, 1, k^{-1}\}$ in case  $g_{1} g_{2}\cdots g_{2n+1} = h^{-1}$ being so $hk \in [g].$ If $hk=1$ clearly $\{\mathcal{E}_{h}, \mathcal{E}_{1}, \mathcal{E}_{h^{-1}}\}\subset \mathcal{E}_{1,[g]}.$ We have showed $\{\mathcal{E}_{h}, \mathcal{E}_{1},  \mathcal{E}_{k}\}\subset \mathcal{E}_{[g]}.$

Next, let us show that $$\{V_{[g]},V_{[g]},\mathcal{E}_{1,[g]}\}\subset \mathcal{E}_{[g]}.$$
\noindent if  $\{\mathcal{E}_{h}, \mathcal{E}_{k}, \mathcal{E}_{1}\} \neq 0$ for some $h, k \in [g]$ then $hk \in \sum^{0}$ and $hk \in \sum^{1}$ $\cup$ $\{1\}$. From here, if $hk \neq 1$ and $\{g_{1},\cdots, g_{2n+1}\}$ is a
connection from $g$ to $h$
then $\{g_{1},\cdots, g_{2n+1}, k, 1\}$ is a
connection from $g$ to $hk$ in
case $g_{1} g_{2}\cdots g_{2n+1} = h$ and  $\{g_{1},\cdots, g_{2n+1}, k^{-1}, 1\}$ in case  $g_{1} g_{2}\cdots g_{2n+1} = h^{-1}$ being so $hk \in [g].$ If $hk=1$ clearly $\{\mathcal{E}_{h},  \mathcal{E}_{h^{-1}},\mathcal{E}_{1}\}\subset \mathcal{E}_{1,[g]}.$ We have showed $\{\mathcal{E}_{h}, \mathcal{E}_{k}, \mathcal{E}_{1} \}\subset \mathcal{E}_{[g]}.$

Finally, let us show
$$\{V_{[g]}, V_{[g]}, V_{[g]}\} \subset \mathcal{E}_{[g]}.$$
Suppose  $\{\mathcal{E}_{h}, \mathcal{E}_{k}, \mathcal{E}_{l}\} \neq 0$ for some $ h, k, l \in [g]$ being so $hk \in  \sum^{0}$ $\cup$ $\{1\}$ and $hkl \in \sum^{1}$ $\cup$ $\{1\}$. If either $hk = 1$ or $hkl = 1$ then  $\{\mathcal{E}_{h}, \mathcal{E}_{k}, \mathcal{E}_{l}\}=\mathcal{E}_{l} \subset V_{[g]}$ or $\{\mathcal{E}_{h}, \mathcal{E}_{k}, \mathcal{E}_{l}\}\subset \mathcal{E}_{1,[g]}$ respectively. Hence let us consider $hk \in  \sum^{0}$ and $hkl \in \sum^{1}$, and
take a connection $\{g_{1}, \cdots, g_{2n+1}\}$ from $g$ to $h$. We
clearly have $\{g_{1}, \cdots, g_{2n+1}, k, l\}$
is a connection from $g$ to $hkl$ in
case $g_{1}\cdots g_{2n+1} = h$ and $\{g_{1} \cdots g_{2n+1}, k^{-1}, l^{-1}\}$ it
is in
case $g_{1}\cdots g_{2n+1} = h^{-1}$. We have showed $hkl \in [g]$ and so $\{\mathcal{E}_{h}, \mathcal{E}_{k}, \mathcal{E}_{l}\}\subset V_{[g]}$,
which concludes the proof.
\epf

\bdefn
With the above notation, we call   $\mathcal{E}_{[g]}$ the \textbf{subtriple} of $\mathcal{E}$ associated to $[g]$.
\edefn

\section{Decompositions}
 We begin this section by showing that for any $g \in \sum^{1}$, the subtriple $I_{[g]}$ is actually
an ideal of $\mathcal{E}$. We need to state some preliminary results.
 \blem \label{lemma 4.1}
 The following assertions hold.

$\rm 1$. If $g, h \in \sum^{1}$  with $gh \in  \pm  \sum^{0}$ $\cup$ $\{1\}$, then $h \in [g]$.

$\rm 2$. If $g, h \in  \sum^{1}$  and $g \in \pm  \sum^{0}$ with $gh \in \pm  \sum^{1}$ $\cup$ $\{1\}$, then $h \in [g]$.

$\rm 3$. If $g, h \in  \sum^{1}$  and $g,h \in \pm  \sum^{0}$ with $gh \in \pm  \sum^{0}$ $\cup$ $\{1\}$, then $h \in [g]$.

$\rm 4$. If $g, \overline{h} \in  \sum^{1}$ such that $\overline{h} \not \in [g]$, then $[\mathcal{E}_{g}, \mathcal{E}_{\overline{h}}] = [L^{0}_{g}, \mathcal{E}_{\overline{h}}] = [L^{0}_{g}
,L^{0}_{\overline{h}}] = 0$.

\bpf
1. If $gh$ = 1, then $h = g^{-1}$ and so $h \sim g$. Suppose $gh \neq 1$. Since $gh \in  \pm  \sum^{0}$,
we have $\{g, h, g^{-1}\}$ is a connection from $g$ to $h$.

2. We can argue similarly with the connection $\{g, 1, (gh)^{-1}\}$.

3. If $gh$ = 1, then $h = g^{-1}$ and so $h \sim g$. Suppose $gh \neq 1$. Since $gh \in  \pm  \sum^{0}$,
we have $\{g, h, g^{-1}\}$ is a connection from $g$ to $h$.

4. Consequence of 1, 2 and 3. \epf

\elem

\blem \label{lemma 4.2}
If $g, \overline{h} \in \sum^{1}$  are not connected, then $\{\mathcal{E}_{g}, \mathcal{E}_{g^{-1}}, \mathcal{E}_{\overline{h}}\} = 0$.
\bpf
If $[\mathcal{E}_{g}, \mathcal{E}_{g^{-1}}] = 0$ it is clear. Suppose then  $[\mathcal{E}_{g}, \mathcal{E}_{g^{-1}}] \neq 0$ and
$\{\mathcal{E}_{g}, \mathcal{E}_{g^{-1}}, \mathcal{E}_{\overline{h}}\} \neq 0$. By Leibniz identity, one gets
$$\{\mathcal{E}_{g}, \mathcal{E}_{g^{-1}}, \mathcal{E}_{\overline{h}}\}=[[\mathcal{E}_{g}, \mathcal{E}_{g^{-1}}], \mathcal{E}_{\overline{h}}]\subset
[\mathcal{E}_{g}, [\mathcal{E}_{g^{-1}}, \mathcal{E}_{\overline{h}}]]+[[\mathcal{E}_{g}, \mathcal{E}_{\overline{h}}], \mathcal{E}_{g^{-1}}].$$
So either $[\mathcal{E}_{g}, [\mathcal{E}_{g^{-1}}, \mathcal{E}_{\overline{h}}]]\neq 0$ or $[[\mathcal{E}_{g}, \mathcal{E}_{\overline{h}}], \mathcal{E}_{g^{-1}}] \neq 0$, contradicting Lemma \ref{lemma 4.1}-4. From here, $\{\mathcal{E}_{g}, \mathcal{E}_{g^{-1}}, \mathcal{E}_{\overline{h}}\} = 0$.
\epf
\elem

\blem \label{lemma 4.3}
For any $g_{0} \in  \sum^{1}$, if $g \in [g_{0}]$ and $h, k \in  \sum^{1}$ $\cup$ $\{1\}$, the following
assertions hold.

$\rm 1$. If $\{\mathcal{E}_{g}, \mathcal{E}_{h}, \mathcal{E}_{k}\} \neq 0$, then $h, k, ghk \in [g_{0}] \cup \{1\}$.

$\rm 2$. If $\{\mathcal{E}_{h},\mathcal{E}_{g} , \mathcal{E}_{k}\} \neq 0$, then $h, k, hgk \in [g_{0}] \cup \{1\}$.

$\rm 3$. If $\{ \mathcal{E}_{h}, \mathcal{E}_{k}, \mathcal{E}_{g}\} \neq 0$, then $h, k, hkg \in [g_{0}] \cup \{1\}$.

\bpf
1. It is easy to see that $[\mathcal{E}_{g}, \mathcal{E}_{h}]\neq 0$ for $g \in [g_{0}]$ and $h \in  \sum^{1}$ $\cup$ $\{1\}$. By Lemma \ref{lemma 4.1}-1, one gets $h\sim g$ in the case $h\neq 1$. From here, $h \in [g_{0}]\cup \{1\}.$ To show $k, ghk \in [g_{0}]\cup \{1\}.$ We distinguish two cases.

Case 1. Suppose $ghk=1$. It is clear that $ghk \in [g_{0}]\cup \{1\}.$ If we have $k\neq 1$ then $gh \in  \sum^{0}$. As $(gh)^{-1}=k$, then $\{g,h,1\}$ would be a connection from $g$ to $k$ and we conclude that $k \in [g_{0}] \cup \{1\}$.

Case 2. Suppose $ghk\neq 1$. We treat separately two cases. If $gh\neq 1$, then $gh \in  \sum^{0}$ and so $\{g, h, k\}$ is a connection from $g$ to $ghk$. Hence $ghk \in [g_{0}]$. In the case $k\neq 1$, we have $\{g, h, (ghk)^{-1}\}$ is a connection from $g$ to $k$. So $k \in [g_{0}]$. Finally, if  $gh=1$, then necessarily $k \in [g_{0}]$. Indeed, if $k$ is not connected to $g$, we would have by Lemma \ref{lemma 4.2} that $\{\mathcal{E}_{g}, \mathcal{E}_{h}, \mathcal{E}_{k}\}=\{\mathcal{E}_{g}, \mathcal{E}_{g^{-1}}, \mathcal{E}_{k}\} =0$, a contradiction. From here $ghk=k \in  [g_{0}]$.

2. This can be proved completely analogously to 1.

3. By Leibniz identity, $\{\mathcal{E}_{h}, \mathcal{E}_{k}, \mathcal{E}_{g}\}=[[\mathcal{E}_{h}, \mathcal{E}_{k}], \mathcal{E}_{g}]\subset [\mathcal{E}_{h}, [\mathcal{E}_{k}, \mathcal{E}_{g}]]+[[\mathcal{E}_{h},\mathcal{E}_{g}], \mathcal{E}_{k}]$. From $\{\mathcal{E}_{h}, \mathcal{E}_{k}, \mathcal{E}_{g}\}\neq 0$, we obtain either $[\mathcal{E}_{h}, [\mathcal{E}_{k}, \mathcal{E}_{g}]]\neq 0$ or $[[\mathcal{E}_{h},\mathcal{E}_{g}], \mathcal{E}_{k}]\neq 0.$ We treat separately two cases.

Case 1. Suppose $[\mathcal{E}_{h}, [\mathcal{E}_{k}, \mathcal{E}_{g}]]\neq 0$, we will show $h, k, hkg \in  [g_{0}]\cup \{1\}.$ First to show $k \in [g_{0}]\cup \{1\}.$ The fact that $[\mathcal{E}_{k}, \mathcal{E}_{g}]\neq 0$ implies by Lemma \ref{lemma 4.1}-1 that $k\sim g$ in the case $k\neq 1$. From here, $k \in [g_{0}]\cup \{1\}.$ Next to show $h\in [g_{0}]\cup \{1\}.$ Indeed, if $h\neq 1$ and suppose $g$ is not connected with $h$, then $h$ is not connected with $k$ in  the case $k\neq 1$. By Lemma \ref{lemma 4.1}-1,   $[\mathcal{E}_{h}, \mathcal{E}_{k}]= 0$ whever $k\neq 1$, contradicting  $[\mathcal{E}_{h}, [\mathcal{E}_{k}, \mathcal{E}_{g}]]\neq 0$.

Next to show if $h \neq 1$ and in the case $k=1$, we also get $h \in [g_{0}]$. Indeed, suppose $g$ is not connected with $h$, in the case $k=1$, one has $\{\mathcal{E}_{h}, \mathcal{E}_{k}, \mathcal{E}_{g}\}=\{\mathcal{E}_{h}, \mathcal{E}_{1}, \mathcal{E}_{g}\}=[[\mathcal{E}_{h}, \mathcal{E}_{1}], \mathcal{E}_{g}]$. From $[\mathcal{E}_{h}, \mathcal{E}_{1}]\subset L_{h}^{0}$ and Lemma \ref{lemma 4.1}-4, one gets $\{\mathcal{E}_{h}, \mathcal{E}_{1}, \mathcal{E}_{g}\}=0$, a contradiction.

Finally, to show $hkg \in  [g_{0}]\cup \{1\}.$ Suppose $hkg=1$ and so $hkg \in  [g_{0}]\cup \{1\}.$ Suppose $hkg\neq 1$, by $[\mathcal{E}_{k}, \mathcal{E}_{g}]\neq 0$, $kg \in  \sum^{1}$ $\cup$ $\{1\}$. If $kg\neq 1$, then $kg \in \sum^{1}$ and so $\{g,k,h\}$ is a connection from $g$ to $gkh$, hence $gkh \in  [g_{0}]$. If $kg=1$, then necessarily $h \in [g_{0}]\cup \{1\}.$ Indeed, if $h\neq 1$ and $g$ is not connected with $h$, by Lemma \ref{lemma 4.1}-4, $\{\mathcal{E}_{h}, \mathcal{E}_{k}, \mathcal{E}_{g}\}=\{\mathcal{E}_{h}, \mathcal{E}_{g^{-1}}, \mathcal{E}_{g}\}=[[\mathcal{E}_{h}, \mathcal{E}_{g^{-1}}], \mathcal{E}_{g}]=0$, contradicting  $\{\mathcal{E}_{h}, \mathcal{E}_{k}, \mathcal{E}_{g}\}\neq 0$. Therefore, we also have $hkg=h \in [g_{0}] \cup \{1\}$.

Case 2. If $[[\mathcal{E}_{h},\mathcal{E}_{g}], \mathcal{E}_{k}]\neq 0,$ we will show $h, k, hgk \in  [g_{0}]\cup \{1\}.$ First to show $h \in
[g_{0}]\cup \{1\}.$ The fact that  $[\mathcal{E}_{h}, \mathcal{E}_{g}]\neq 0$ implies by Lemma \ref{lemma 4.1}-1 that $h\sim g$ in the case $h\neq 1$. From here, $h \in [g_{0}]\cup \{1\}.$

Next to show $k \in [g_{0}]\cup \{1\}.$ Indeed, if $k\neq 1$ and $g$ is not connected with $k$, then $h$ is not connected with $k$ in the case $h\neq 1$. By Lemma \ref{lemma 4.1}-1,  $[\mathcal{E}_{h}, \mathcal{E}_{k}]=0$ whenever $h\neq 1$, contradicting $\{\mathcal{E}_{h}, \mathcal{E}_{k}, \mathcal{E}_{g}\}\neq 0$. Similarly, it is easy to show if $k\neq 1$ and in the case $h=1$, we can obtain $k \in [g_{0}]$.

Finally, to show $hkg \in [g_{0}]\cup \{1\}.$ Suppose $hkg=1$ and so $hkg \in [g_{0}]\cup \{1\}.$ Suppose  $hkg\neq 1$, by $ [\mathcal{E}_{h}, \mathcal{E}_{g}]\neq 0$, one has $hg \in  \sum^{0} \cup \{1\}.$ If $hg\neq 1$, then $hg \in  \sum^{0} $ and so $\{g,h,k\}$ is a connection from $g$ to $hkg$. Hence $hkg \in [g_{0}]$. If $hg=1$, then necessarily $k \in [g_{0}]\cup \{1\}.$ Indeed, if $k\neq 1$ and $g$ is not connected with $k$, by Lemma \ref{lemma 4.1}-4, $\{\mathcal{E}_{h}, \mathcal{E}_{k}, \mathcal{E}_{g}\}=\{\mathcal{E}_{g^{-1}}, \mathcal{E}_{k}, \mathcal{E}_{g}\}=[[\mathcal{E}_{g^{-1}}, \mathcal{E}_{k}], \mathcal{E}_{g}]=0$, contradicting  $\{\mathcal{E}_{h}, \mathcal{E}_{k}, \mathcal{E}_{g}\}\neq 0$. Therefore, we also have $hkg=k \in [g_{0}]\cup \{1\}.$
\epf

\elem

\blem \label{lemma 4.4}
For any $g_{0} \in  \sum^{1}$, if $g$, $k \in [g_{0}]$,  $h \in  [g_{0}]$  $\cup$ $\{1\}$ with $ghk=1$ and $l,m \in \sum^{1}$ $\cup$ $\{1\}$, the following assertions hold.

$\rm 1$. If $\{\{\mathcal{E}_{g}, \mathcal{E}_{h}, \mathcal{E}_{k}\}, \mathcal{E}_{l},\mathcal{E}_{m}\}\neq 0$, then $l, m, lm \in [g_{0}] \cup \{1\}$.

$\rm 2$. If $\{\mathcal{E}_{l},\{\mathcal{E}_{g},\mathcal{E}_{h} , \mathcal{E}_{k}\},\mathcal{E}_{m}\} \neq 0$, then $l, m, lm \in [g_{0}] \cup \{1\}$.

$\rm 3$.  If $\{\mathcal{E}_{l},\mathcal{E}_{m},\{\mathcal{E}_{g}, \mathcal{E}_{h}, \mathcal{E}_{k}\}\} \neq 0$, then $l, m, lm \in [g_{0}] \cup \{1\}$.

\bpf
1. By (\ref{VIP2}), one gets
\begin{align*}
0\neq  \{\{\mathcal{E}_{g}, \mathcal{E}_{h}, \mathcal{E}_{k}\}, \mathcal{E}_{l},\mathcal{E}_{m}\} &\subset  \{\mathcal{E}_{g}, \mathcal{E}_{h}, \{\mathcal{E}_{k}, \mathcal{E}_{l},\mathcal{E}_{m}\}\}+\{\{\mathcal{E}_{g}, \mathcal{E}_{h}, \mathcal{E}_{l}\}, \mathcal{E}_{k},\mathcal{E}_{m}\}\\
&+\{\{\mathcal{E}_{g}, \mathcal{E}_{h}, \mathcal{E}_{m}\}, \mathcal{E}_{k},\mathcal{E}_{l}\}+\{\{\mathcal{E}_{g}, \mathcal{E}_{h}, \mathcal{E}_{m}\}, \mathcal{E}_{l},\mathcal{E}_{k}\},
\end{align*}
\noindent any of the above four summands is nonzero. In order to prove $l$, $m$, $lm \in   [g_{0}]$  $\cup$ $\{1\}$, we will consider four cases.

Case 1. Suppose $\{\mathcal{E}_{g}, \mathcal{E}_{h}, \{\mathcal{E}_{k}, \mathcal{E}_{l},\mathcal{E}_{m}\}\}\neq 0.$
 As $k \in [g_{0}]$ and $\{\mathcal{E}_{k}, \mathcal{E}_{l},\mathcal{E}_{m}\}\neq 0$, Lemma \ref{lemma 4.3}-1 shows that $l, m, lm $ are
connected with $k$ in the case of being nonzero roots and so $l, m, klm \in   [g_{0}]$  $\cup$ $\{1\}$.
If $klm=1$, then $lm=k^{-1}\in [g_{0}]$. If $klm \neq 1$, taking into account $0\neq
\{\mathcal{E}_{g}, \mathcal{E}_{h}, \{\mathcal{E}_{k}, \mathcal{E}_{l},\mathcal{E}_{m}\}\}\subset \{\mathcal{E}_{g}, \mathcal{E}_{h}, \mathcal{E}_{klm}\}$, Lemma \ref{lemma 4.3}-1 gives
 us that $ghklm=lm \in  [g_{0}]$.

Case 2. Suppose $\{\{\mathcal{E}_{g}, \mathcal{E}_{h}, \mathcal{E}_{l}\}, \mathcal{E}_{k},\mathcal{E}_{m}\}\neq 0$. It is clear that  $\{\mathcal{E}_{g}, \mathcal{E}_{h}, \mathcal{E}_{l}\}\neq 0$.  As $g \in [g_{0}]$, by Lemma \ref{lemma 4.3}-1, one gets $l \in  [g_{0}]$  $\cup$ $\{1\}$. It is obvious that $0\neq \{\{\mathcal{E}_{g}, \mathcal{E}_{h}, \mathcal{E}_{l}\}, \mathcal{E}_{k},\mathcal{E}_{m}\}\subset \{\mathcal{E}_{ghl}, \mathcal{E}_{k},\mathcal{E}_{m}\}$. As $k \in [g_{0}]$,  by Lemma \ref{lemma 4.3}-2, one gets $m \in  [g_{0}]$  $\cup$ $\{1\}$ and $ghlkm=lm \in  [g_{0}]$  $\cup$ $\{1\}$.

Case 3. Suppose $\{\{\mathcal{E}_{g}, \mathcal{E}_{h}, \mathcal{E}_{m}\}, \mathcal{E}_{k},\mathcal{E}_{l}\}\neq 0$. It is easy to see that $\{\mathcal{E}_{g}, \mathcal{E}_{h}, \mathcal{E}_{m}\}\neq 0$. As $g \in [g_{0}]$, by Lemma \ref{lemma 4.3}-1, one gets $m \in  [g_{0}]$  $\cup$
$\{1\}$. Note that $0\neq \{\{\mathcal{E}_{g}, \mathcal{E}_{h}, \mathcal{E}_{m}\}, \mathcal{E}_{k},\mathcal{E}_{l}\}\subset \{\mathcal{E}_{ghm}, \mathcal{E}_{k}, \mathcal{E}_{l}\}$. As  $k \in [g_{0}]$,  by Lemma \ref{lemma 4.3}-2, one gets $l \in  [g_{0}]$  $\cup$ $\{1\}$ and $ghklm=lm \in  [g_{0}]$  $\cup$ $\{1\}$.

Case 4. Suppose $\{\{\mathcal{E}_{g}, \mathcal{E}_{h}, \mathcal{E}_{m}\}, \mathcal{E}_{l},\mathcal{E}_{k}\}\neq 0$. It is clear that $\{\mathcal{E}_{g}, \mathcal{E}_{h}, \mathcal{E}_{m}\}\neq 0$. As $g \in [g_{0}]$, by Lemma \ref{lemma 4.3}-1, one gets $m \in  [g_{0}]$  $\cup$ $\{1\}$. Note that
$0\neq \{\{\mathcal{E}_{g}, \mathcal{E}_{h}, \mathcal{E}_{m}\}, \mathcal{E}_{l},\mathcal{E}_{k}\}\subset \{\mathcal{E}_{ghm},  \mathcal{E}_{l},\mathcal{E}_{k}\}$. As  $k \in [g_{0}]$,  by Lemma \ref{lemma 4.3}-3, one gets $l \in  [g_{0}]$  $\cup$ $\{1\}$ and $ghmlk=lm \in  [g_{0}]$  $\cup$ $\{1\}$.

2. By Proposition \ref{38888} (3), we obtain that
\begin{align*}
0&\neq  \{\mathcal{E}_{l}, \{\mathcal{E}_{g}, \mathcal{E}_{h}, \mathcal{E}_{k}\},\mathcal{E}_{m}\}\\
 &\subset  \{\{\mathcal{E}_{l}, \mathcal{E}_{k}, \mathcal{E}_{m}\}, \mathcal{E}_{h},\mathcal{E}_{g}\}+\{\{\mathcal{E}_{l}, \mathcal{E}_{k}, \mathcal{E}_{m}\}, \mathcal{E}_{g},\mathcal{E}_{h}\}\\
&+\{\{\mathcal{E}_{l}, \mathcal{E}_{h}, \mathcal{E}_{g}\}, \mathcal{E}_{k},\mathcal{E}_{m}\}+\{\{\mathcal{E}_{l}, \mathcal{E}_{g}, \mathcal{E}_{h}\}, \mathcal{E}_{k},\mathcal{E}_{m}\}+\{\mathcal{E}_{l}, \mathcal{E}_{k}, \{\mathcal{E}_{g}, \mathcal{E}_{h},\mathcal{E}_{m}\}\},
\end{align*}
\noindent any of the above five summands is nonzero. Suppose $\{\{\mathcal{E}_{l}, \mathcal{E}_{k}, \mathcal{E}_{m}\}, \mathcal{E}_{h},\mathcal{E}_{g}\}\neq 0$, it is obvious $\{\mathcal{E}_{l}, \mathcal{E}_{k}, \mathcal{E}_{m}\}\neq 0$. As $k \in [g_{0}]$,  by Lemma \ref{lemma 4.3}-2, one gets $l,lkm \in  [g_{0}]$  $\cup$ $\{1\}$. Note that $0\neq \{\{\mathcal{E}_{l}, \mathcal{E}_{k}, \mathcal{E}_{m}\}, \mathcal{E}_{h},\mathcal{E}_{g}\}\subset \{\mathcal{E}_{lkm},\mathcal{E}_{h},\mathcal{E}_{g}\}$. As  $g \in [g_{0}]$, by Lemma \ref{lemma 4.3}-3, one gets $lkmhg=lm \in  [g_{0}]$  $\cup$ $\{1\}$. If $\{\{\mathcal{E}_{l}, \mathcal{E}_{k}, \mathcal{E}_{m}\}, \mathcal{E}_{g},\mathcal{E}_{h}\}\neq 0$,
$\{\{\mathcal{E}_{l}, \mathcal{E}_{h}, \mathcal{E}_{g}\}, \mathcal{E}_{k},\mathcal{E}_{m}\}\neq 0$, $\{\{\mathcal{E}_{l}, \mathcal{E}_{g}, \mathcal{E}_{h}\}, \mathcal{E}_{k},\mathcal{E}_{m}\}\neq 0$, $\{\mathcal{E}_{l}, \mathcal{E}_{k}, \{\mathcal{E}_{g}, \mathcal{E}_{h},\mathcal{E}_{m}\}\}\neq 0$, a similar argument gives us $l,m,lm \in [g_{0}]$  $\cup$ $\{1\}$.

3. By Proposition \ref{38888} (3), we obtain that
\begin{align*}
0&\neq  \{\mathcal{E}_{l}, \mathcal{E}_{m}, \{\mathcal{E}_{g}, \mathcal{E}_{h},\mathcal{E}_{k}\}\}\\
 &\subset  \{\{\mathcal{E}_{l}, \mathcal{E}_{m}, \mathcal{E}_{k}\}, \mathcal{E}_{h},\mathcal{E}_{g}\}+\{\{\mathcal{E}_{l}, \mathcal{E}_{m}, \mathcal{E}_{k}\}, \mathcal{E}_{g},\mathcal{E}_{h}\}\\
&+\{\{\mathcal{E}_{l}, \mathcal{E}_{h}, \mathcal{E}_{g}\}, \mathcal{E}_{m},\mathcal{E}_{k}\}+\{\{\mathcal{E}_{l}, \mathcal{E}_{g}, \mathcal{E}_{h}\}, \mathcal{E}_{m},\mathcal{E}_{k}\}+\{\mathcal{E}_{l}, \{\mathcal{E}_{g}, \mathcal{E}_{h}, \mathcal{E}_{m}\},\mathcal{E}_{k}\},
\end{align*}
\noindent any of the above five summands is nonzero. Suppose $\{\{\mathcal{E}_{l}, \mathcal{E}_{m}, \mathcal{E}_{k}\}, \mathcal{E}_{h},\mathcal{E}_{g}\}\neq 0$, one easily gets $\{\mathcal{E}_{l}, \mathcal{E}_{m}, \mathcal{E}_{k}\}\neq 0$. As $k \in [g_{0}]$,  by Lemma \ref{lemma 4.3}-3, one gets $l,lmk \in  [g_{0}]$  $\cup$ $\{1\}$. Note that $0\neq  [[\mathcal{E}_{l}, \mathcal{E}_{m}, \mathcal{E}_{k}], \mathcal{E}_{h},\mathcal{E}_{g}]\subset [\mathcal{E}_{lmk},\mathcal{E}_{h},\mathcal{E}_{g}]$. As $g \in [g_{0}]$,  by Lemma \ref{lemma 4.3}-3, one gets $lmkhg=lm \in  [g_{0}]$  $\cup$ $\{1\}$. If $\{\{\mathcal{E}_{l}, \mathcal{E}_{m}, \mathcal{E}_{k}\}, \mathcal{E}_{g},\mathcal{E}_{h}\}\neq 0$,
$\{\{\mathcal{E}_{l}, \mathcal{E}_{h}, \mathcal{E}_{g}\}, \mathcal{E}_{m},\mathcal{E}_{k}\}\neq 0$, $\{\{\mathcal{E}_{l}, \mathcal{E}_{g}, \mathcal{E}_{h}\}, \mathcal{E}_{m},\mathcal{E}_{k}\}\neq 0$ or $\{\mathcal{E}_{l}, \{\mathcal{E}_{g}, \mathcal{E}_{h}, \mathcal{E}_{m}\},\mathcal{E}_{k}\}\neq 0$, a similar argument gives us $l,m,lm \in [g_{0}]$  $\cup$ $\{1\}$.
\epf
\elem

\blem \label{lemma 4.5}
For any $g_{0} \in  \sum^{1}$, if $g$, $k \in [g_{0}]$,  $h \in  [g_{0}]$  $\cup$ $\{1\}$ with $ghk=1$ and $\overline{h} \not \in  [g_{0}]$, the following assertions hold.

$\rm 1$.  $[\{\mathcal{E}_{g}, \mathcal{E}_{h}, \mathcal{E}_{k}\}, \mathcal{E}_{\overline{h}}]=0$.

$\rm 2$. $[\{\mathcal{E}_{g}, \mathcal{E}_{h}, \mathcal{E}_{k}\}, L_{\overline{h}}^{0}]=0$.

$\rm 3$.  $\{\{\mathcal{E}_{g}, \mathcal{E}_{h}, \mathcal{E}_{k}\}, \mathcal{E}_{1}, \mathcal{E}_{\overline{h}}\}=0$.
\elem
\bpf
1. By Leibniz identity, we have
\begin{align}\label{cao4.5}
 [\{\mathcal{E}_{g}, \mathcal{E}_{h}, \mathcal{E}_{k}\}, \mathcal{E}_{\overline{h}}]=[[[\mathcal{E}_{g}, \mathcal{E}_{h}], \mathcal{E}_{k}], \mathcal{E}_{\overline{h}}]\subset  [[\mathcal{E}_{g}, \mathcal{E}_{h}],[\mathcal{E}_{k}, \mathcal{E}_{\overline{h}}]]+[[[\mathcal{E}_{g}, \mathcal{E}_{h}], \mathcal{E}_{\overline{h}}],\mathcal{E}_{k}].
\end{align}

\noindent Let us consider the first summand in (\ref{cao4.5}). As  $k \in [g_{0}]$, for $\overline{h} \not \in  [g_{0}]$, by Lemma \ref{lemma 4.1}-4, one gets $[\mathcal{E}_{k}, \mathcal{E}_{\overline{h}}]=0$. Therefore  $[[\mathcal{E}_{g}, \mathcal{E}_{h}],[\mathcal{E}_{k}, \mathcal{E}_{\overline{h}}]]=0$. Let us now consider the second summand in  (\ref{cao4.5}), it is sufficient to verify that $$[[[\mathcal{E}_{g}, \mathcal{E}_{h}], \mathcal{E}_{\overline{h}}],\mathcal{E}_{k}]=0.$$
\noindent To do so, we first assert that $[[\mathcal{E}_{g}, \mathcal{E}_{h}], \mathcal{E}_{\overline{h}}]=0.$ Indeed, by  Leibniz identity, we have \begin{align}\label{cao4.6}
 [[\mathcal{E}_{g}, \mathcal{E}_{h}],  \mathcal{E}_{\overline{h}}]\subset  [\mathcal{E}_{g}, [\mathcal{E}_{h},\mathcal{E}_{\overline{h}}]]+[[\mathcal{E}_{g}, \mathcal{E}_{\overline{h}}], \mathcal{E}_{h}],
\end{align}
\noindent where  $g \in [g_{0}]$,  $h \in  [g_{0}]$  $\cup$ $\{1\},$ $\overline{h} \not \in  [g_{0}]$.  Let us consider the first summand in (\ref{cao4.6}), if $h \not \in \{1\}$, then $h \in  [g_{0}]$. By Lemma \ref{lemma 4.1}-1, one gets $[\mathcal{E}_{h},\mathcal{E}_{\overline{h}}]=0$. That is, $[\mathcal{E}_{g}, [\mathcal{E}_{h}, \mathcal{E}_{\overline{h}}]]=0$. If $h=1$, we have $[\mathcal{E}_{h},\mathcal{E}_{\overline{h}}]\subset L_{\overline{h}}^{0}$. By Lemma  \ref{lemma 4.1}-4, one gets $[\mathcal{E}_{g}, [\mathcal{E}_{h},\mathcal{E}_{\overline{h}}]]\subset [\mathcal{E}_{g},  L_{\overline{h}}^{0}]=0$. Therefore $[\mathcal{E}_{g}, [\mathcal{E}_{h},\mathcal{E}_{\overline{h}}]]=0$. Let us consider the second summand in (\ref{cao4.6}), it is sufficient to verify that
$[[\mathcal{E}_{g}, \mathcal{E}_{\overline{h}}], \mathcal{E}_{h}]=0$. Indeed, for  $g \in [g_{0}]$
 and $\overline{h} \not \in  [g_{0}]$. By Lemma  \ref{lemma 4.1}-4, $[[\mathcal{E}_{g}, \mathcal{E}_{\overline{h}}], \mathcal{E}_{h}]=0$.

2. Note that
\begin{align}\label{cao4.7}
[\{\mathcal{E}_{g}, \mathcal{E}_{h}, \mathcal{E}_{k}\}, L_{\overline{h}}^{0}] \subset [[\mathcal{E}_{g}, \mathcal{E}_{h}], [\mathcal{E}_{k}, L_{\overline{h}}^{0}]]+[[[\mathcal{E}_{g}, \mathcal{E}_{h}],  L_{\overline{h}}^{0}], \mathcal{E}_{k}].
\end{align}

\noindent Let us consider the first summand in (\ref{cao4.7}). As $k\neq 1$, one gets  $[[\mathcal{E}_{g}, \mathcal{E}_{h}], [\mathcal{E}_{k}, L_{\overline{h}}^{0}]]=0$ by Lemma  \ref{lemma 4.1}-4. Let us consider the second summand in (\ref{cao4.7}). By Leibniz identity, we obtain $[[\mathcal{E}_{g}, \mathcal{E}_{h}],  L_{\overline{h}}^{0}]=0$, so $[[[\mathcal{E}_{g}, \mathcal{E}_{h}],  L_{\overline{h}}^{0}], \mathcal{E}_{k}]=0$.

3. It is a consequence of Lemma  \ref{lemma 4.5}-1, 2 and
$$\{\{\mathcal{E}_{g}, \mathcal{E}_{h}, \mathcal{E}_{k}\},
\mathcal{E}_{1}, \mathcal{E}_{\overline{h}}\}\subset
[\{\mathcal{E}_{g}, \mathcal{E}_{h}, \mathcal{E}_{k}\},
[\mathcal{E}_{1}, \mathcal{E}_{\overline{h}}]]+[[\{\mathcal{E}_{g},
\mathcal{E}_{h}, \mathcal{E}_{k}\},
\mathcal{E}_{\overline{h}}],\mathcal{E}_{1}].$$ Thus the lemma
follows.\epf

\bdefn\label{ 3.16789}
A Leibniz triple system $\mathcal{E}$ is said to be \textbf{simple} if its product is nonzreo and its only ideals are $\{0\}$, $J$ and $\mathcal{E}$.
\edefn

It should be noted that the above definition agrees with the definition of a simple Lie triple system, since $J=\{0\}$ in this case.

\bthm\label{theorem 1}
The following assretions hold.

$\rm 1$. For any $g_{0} \in \sum^{1}$, the subtriple $\mathcal{E}_{[g_{0}]}=\mathcal{E}_{1,[g_{0}]}\oplus V_{[g_{0}]}$ of $\mathcal{E}$ associated to $[g_{0}]$ is an ideal of $\mathcal{E}$.

$\rm 2$.  If  $\mathcal{E}$ is simple, then $\sum^{1}$ has all of its elements connected and
$$\mathcal{E}_{1}=\sum_{g \in \sum^{1}, \ h \in \sum^{1}\cup \{1\}}\{\mathcal{E}_{g},\mathcal{E}_{h},\mathcal{E}_{(gh)^{-1}}\}.$$
\ethm

\bpf 1.
Recall that
\begin{align}\label{cao 5}
\mathcal{E}_{1,[g_{0}]} := span_{\mathbb{K}}\{\{\mathcal{E}_{g}, \mathcal{E}_{h}, \mathcal{E}_{(gh)^{-1}}\} : g \in [g_{0}], h \in [g_{0}] \cup\{1\}\}   \subset \mathcal{E}_{1}.
\end{align}
\noindent In order to complete the proof, it is sufficient to show that
$$\{\mathcal{E}_{[g_{0}]}, \mathcal{E}, \mathcal{E}\}+\{\mathcal{E}, \mathcal{E}_{[g_{0}]}, \mathcal{E}\}+\{\mathcal{E}, \mathcal{E}, \mathcal{E}_{[g_{0}]}\}\subset \mathcal{E}_{[g_{0}]}.$$
\noindent We first check that $\{\mathcal{E}_{[g_{0}]}, \mathcal{E}, \mathcal{E}\}\subset \mathcal{E}_{[g_{0}]}.$ We easily have $\{\mathcal{E}_{1,[g_{0}]}, \mathcal{E}, \mathcal{E}\}\subset \mathcal{E}_{[g_{0}]}$. (\ref{cao 5}) together with Lemma \ref{lemma 4.4} imply
$$\{\mathcal{E}_{1,[g_{0}]}, \mathcal{E}_{1},  \mathcal{E}_{g}\}+\{\mathcal{E}_{1,[g_{0}]}, \mathcal{E}_{g}, \mathcal{E}_{1}\}+\{\mathcal{E}_{1,[g_{0}]}, \mathcal{E}_{g}, \mathcal{E}_{h}\}\subset \mathcal{E}_{[g_{0}]},$$
\noindent for any $g,h \in \sum^{1}$. From here,
\begin{align}\label{cao 6}
\{\mathcal{E}_{1,[g_{0}]}, \mathcal{E},  \mathcal{E}\}= \{\mathcal{E}_{1,[g_{0}]}, \mathcal{E}_{1}\oplus(\oplus_{g \in \sum^{1}}\mathcal{E}_{g}), \mathcal{E}_{1}\oplus(\oplus_{h \in \sum^{1}}\mathcal{E}_{h}) \}\subset \mathcal{E}_{[g_{0}]}.
\end{align}
\noindent Since $V_{[g_{0}]}:=\oplus_{g \in [g_{0}]}\mathcal{E}_{g}$, we have by Lemma \ref{lemma 4.3} and (\ref{cao 5}) that
$$\{\oplus_{g \in [g_{0}]}\mathcal{E}_{g}, \mathcal{E}_{1}, \mathcal{E}_{1}\}+\{\oplus_{g \in [g_{0}]}\mathcal{E}_{g}, \mathcal{E}_{1}, \mathcal{E}_{h}\}+\{\oplus_{g \in [g_{0}]}\mathcal{E}_{g}, \mathcal{E}_{h}, \mathcal{E}_{1} \}+\{\oplus_{g \in [g_{0}]}\mathcal{E}_{g}, \mathcal{E}_{h}, \mathcal{E}_{k} \}\subset \mathcal{E}_{g_{0}},$$
\noindent for any $h,k \in \sum^{1}$. So
\begin{align}\label{cao 7}
\{V_{[g_{0}]}, \mathcal{E},  \mathcal{E}\}= \{\oplus_{g \in [g_{0}]}\mathcal{E}_{g}, \mathcal{E}_{1}\oplus(\oplus_{h \in \sum^{1}}\mathcal{E}_{h}), \mathcal{E}_{1}\oplus(\oplus_{k \in \sum^{1}}\mathcal{E}_{k}) \}\subset \mathcal{E}_{[g_{0}]}.
\end{align}
\noindent From (\ref{cao 6}) and (\ref{cao 7}), we have
$$\{\mathcal{E}_{[g_{0}]}, \mathcal{E}, \mathcal{E}\}=\{\mathcal{E}_{1,[g_{0}]}\oplus V_{[g_{0}]}, \mathcal{E}, \mathcal{E}\}\subset \mathcal{E}_{[g_{0}]},$$
\noindent and so    $\mathcal{E}_{[g_{0}]}$ is an ideal of $\mathcal{E}$.

By Lemmas \ref{lemma 4.3} and \ref{lemma 4.4}, a similar argument gives us $[\mathcal{E}, \mathcal{E}_{[g_{0}]}, \mathcal{E}]\subset \mathcal{E}_{[g_{0}]}$ and $[ \mathcal{E}, \mathcal{E}, \mathcal{E}_{[g_{0}]}]\subset \mathcal{E}_{[g_{0}]}.$ Consequently, this proves $\mathcal{E}_{[g_{0}]}$ is an ideal of $\mathcal{E}$.

2. The simplicity of $\mathcal{E}$ implies $\mathcal{E}_{[g_{0}]} \in \{J, \mathcal{E}\}$ for any $g \in [g_{0}]$. If $g \in  [g_{0}]$ is such that $\mathcal{E}_{[g_{0}]}=\mathcal{E}$. Then $ [g_{0}]=\sum^{1}$. Hence, $\mathcal{E}$ has all its nonzero roots connected. Otherwise, if $\mathcal{E}_{[g_{0}]}=J$ for any $g \in [g_{0}]$ then $\mathcal{E}_{[g_{0}]}=\mathcal{E}_{[\alpha_{0}]}$ for any $g_{0}, \alpha_{0} \in \sum^{1}$ and so $[g_{0}]= \sum^{1}$, we also conclude that $\mathcal{E}$ has all its nonzero roots connected.
\epf

\bthm\label{theorem 2}
For a linear complement $U$ of $span_{\mathbb{K}}\{\{\mathcal{E}_{g}, \mathcal{E}_{h}, \mathcal{E}_{(gh)^{-1}}\}: g \in \sum^{1}, h \in \sum^{1}\cup\{1\}\}$ in   $\mathcal{E}_{1}$, we have
$$\mathcal{E}=U+\sum_{[g]\in \sum^{1}/ \sim}I_{[g]},$$
\noindent where any $I_{[g]}$ is one of the ideals described in Theorem \ref{theorem 1}, which also satisfy $[I_{[g]}, \mathcal{E}, I_{[h]}]=[I_{[g]}, I_{[h]}, \mathcal{E}]=[\mathcal{E}, I_{[g]},  I_{[h]}]=0$ if $[g]\neq [h]$.
\ethm

\bpf
By proposition \ref{6787777}, we can consider the quotient set  $\sum^{1}/ \sim:=\{[g]: g \in \sum^{1}\}.$ We have $I_{[g]}$ is well defined and by Theorem \ref{theorem 1}-1 an ideal of $ \mathcal{E}$. Therefore
$$\mathcal{E}=U+\sum_{[g]\in \sum^{1}/ \sim}I_{[g]}.$$
\noindent Next, it is sufficient to show that $\{I_{[g]}, \mathcal{E}, I_{[h]}\}=0$ if  $[g]\neq [h]$. Note that,
\begin{align*}
 \{I_{[g]}, \mathcal{E}, I_{[h]}\}&=\{\mathcal{E}_{1,[g]}\oplus V_{[g]},\mathcal{E}_{1}\oplus(\oplus_{k \in \sum^{1}}\mathcal{E}_{k}),\mathcal{E}_{1,[h]}\oplus V_{[h]}\}\\
&=\{\mathcal{E}_{1,[g]},\mathcal{E}_{1}, \mathcal{E}_{1,[h]}\}+\{\mathcal{E}_{1,[g]},\mathcal{E}_{1}, V_{[h]}\}+\{\mathcal{E}_{1,[g]}, \oplus_{k \in \sum^{1}}\mathcal{E}_{k},\mathcal{E}_{1,[h]} \}\\
&+\{\mathcal{E}_{1,[g]}, \oplus_{k \in \sum^{1}}\mathcal{E}_{k},V_{[h]} \}+\{V_{[g]},\mathcal{E}_{1}, \mathcal{E}_{1,[h]}\}+\{V_{[g]},\mathcal{E}_{1}, V_{[h]}\}\\
&+\{V_{[g]},\oplus_{k \in \sum^{1}}\mathcal{E}_{k},\mathcal{E}_{1,[h]}\}+\{V_{[g]},\oplus_{k \in \sum^{1}}\mathcal{E}_{k},V_{[h]}\}.
\end{align*}
\noindent Hence, if $[\alpha]\neq [\beta]$, by Lemmas \ref{lemma 4.3} and \ref{lemma 4.4}, it is easy to see
\begin{align*}
\{\mathcal{E}_{1,[g]},\mathcal{E}_{1}, V_{[h]}\}&=\{\mathcal{E}_{1,[g]}, \oplus_{k \in \sum^{1}}\mathcal{E}_{k},V_{[h]} \}=\{V_{[g]},\mathcal{E}_{1}, \mathcal{E}_{1,[h]}\}=\{V_{[g]},\mathcal{E}_{1}, V_{[h]}\}\\
&=\{V_{[g]},\oplus_{k \in \sum^{1}}\mathcal{E}_{k},\mathcal{E}_{1,[h]}\}=\{V_{[g]},\oplus_{k \in \sum^{1}}\mathcal{E}_{k},V_{[h]}\}=0.
\end{align*}

Next, we will show that $\{\mathcal{E}_{1,[g]}, \oplus_{k \in \sum^{1}}\mathcal{E}_{k},\mathcal{E}_{1,[h]} \}=0$. Indeed, for $\{\mathcal{E}_{\alpha_{1}},\mathcal{E}_{\alpha_{2}},\mathcal{E}_{(\alpha_{1}\alpha_{2})^{-1}}\}\in \mathcal{E}_{1,[g]}$ with $\alpha_{1} \in [g]$,
$\alpha_{2} \in [g]\cup \{1\}$ and for $\{\mathcal{E}_{\beta_{1}},\mathcal{E}_{\beta_{2}},\mathcal{E}_{(\beta_{1}\beta_{2})^{-1}}\}\in \mathcal{E}_{1,[h]}$ with $\beta_{1} \in [h]$,
$\beta_{2} \in [h]\cup \{1\},$ by Proposition \ref{38888} (3) and Lemma \ref{lemma 4.4}, one gets
\begin{align*}
& \{\{\mathcal{E}_{\alpha_{1}},\mathcal{E}_{\alpha_{2}},\mathcal{E}_{(\alpha_{1}\alpha_{2})^{-1}}\}, \oplus_{k \in \sum^{1}}\mathcal{E}_{k},  \{\mathcal{E}_{\beta_{1}},\mathcal{E}_{\beta_{2}},\mathcal{E}_{(\beta_{1}\beta_{2})^{-1}}\}\}\\
 &\subset \{\{\{\mathcal{E}_{\alpha_{1}},\mathcal{E}_{\alpha_{2}},\mathcal{E}_{(\alpha_{1}\alpha_{2})^{-1}}\}, \oplus_{k \in \sum^{1}}\mathcal{E}_{k},  \mathcal{E}_{(\beta_{1}\beta_{2})^{-1}}\}, \mathcal{E}_{\beta_{2}}, \mathcal{E}_{\beta_{1}}\}\\
&+\{\{\{\mathcal{E}_{\alpha_{1}},\mathcal{E}_{\alpha_{2}},\mathcal{E}_{(\alpha_{1}\alpha_{2})^{-1}}\}, \oplus_{k \in \sum^{1}}\mathcal{E}_{k},  \mathcal{E}_{(\beta_{1}\beta_{2})^{-1}}\}, \mathcal{E}_{\beta_{1}}, \mathcal{E}_{\beta_{2}}\}\\
&+\{\{\{\mathcal{E}_{\alpha_{1}},\mathcal{E}_{\alpha_{2}},\mathcal{E}_{(\alpha_{1}\alpha_{2})^{-1}}\},  \mathcal{E}_{\beta_{2}} \mathcal{E}_{\beta_{1}}\},\oplus_{k \in \sum^{1}}\mathcal{E}_{k},  \mathcal{E}_{(\beta_{1}\beta_{2})^{-1}}\}, \\
&+\{\{\{\mathcal{E}_{\alpha_{1}},\mathcal{E}_{\alpha_{2}},\mathcal{E}_{(\alpha_{1}\alpha_{2})^{-1}}\},  \mathcal{E}_{\beta_{1}} \mathcal{E}_{\beta_{2}}\},\oplus_{k \in \sum^{1}}\mathcal{E}_{k},  \mathcal{E}_{(\beta_{1}\beta_{2})^{-1}}\}, \\
&+\{\{\mathcal{E}_{\alpha_{1}},\mathcal{E}_{\alpha_{2}},\mathcal{E}_{(\alpha_{1}\alpha_{2})^{-1}}\},  \{\mathcal{E}_{\beta_{1}} \mathcal{E}_{\beta_{2}},\oplus_{k \in \sum^{1}}\mathcal{E}_{k}\},  \mathcal{E}_{(\beta_{1}\beta_{2})^{-1}}\}\\
&=0.
\end{align*}

\noindent A similar method gives that $\{\mathcal{E}_{1,[g]},\mathcal{E}_{1}, \mathcal{E}_{1,[h]}\}=0$. So we prove that  $\{I_{[g]}, \mathcal{E}, I_{[h]}\}=0$ if  $[g]\neq [h]$. A similar argument gives that $\{I_{[g]}, I_{[h]}, \mathcal{E}\}=\{\mathcal{E}, I_{[g]},  I_{[h]}\}=0$ if $[g]\neq [h]$.
\epf

\bdefn
The \textbf{annihilator} of a Leibniz triple system $\mathcal{E}$ is the set $\mathrm{Ann}(\mathcal{E})=\{x \in \mathcal{E}: \{x, \mathcal{E}, \mathcal{E}\}+ \{\mathcal{E}, x,  \mathcal{E}\}+\{\mathcal{E},\mathcal{E},x\}=0\}$.
\edefn

\bdefn We will say that $\mathcal{E}_{1}$ is \textbf{tight} if
$\mathcal{E}_{1}= {\rm span}_{\mathbb{K}}\{\{\mathcal{E}_{g},
\mathcal{E}_{h}, \mathcal{E}_{(gh)^{-1}}\} : g \in \sum^{1}, h \in
\sum^{1} \cup\{1\}\}.$ \edefn

\bcor
If $\mathrm{Ann}(\mathcal{E})=0$ and $\mathcal{E}_{1}$ is tight, then $\mathcal{E}$ is the direct sum of the ideals given in Theorem \ref{theorem 1}-1,
$$\mathcal{E}=\oplus_{[g] \in \sum^{1}/\sim}I_{[g]}.$$
\ecor
\bpf
From the fact that  $\mathcal{E}_{1}$ is tight, we  clearly have $$\mathcal{E}=\oplus_{[g] \in \sum^{1}/\sim}I_{[g]}.$$ To finish, we show the direct character of the sum. Given $x \in I_{[g]}\cap (\sum_{[h]\in(\sum^{1}/\sim)\setminus [g]} I_{[h]})$ we have from the fact   $\{I_{[g]}, \mathcal{E}, I_{[h]}\}=0$ that
$$\{x,\mathcal{E}, I_{[g]} \}+\{x,\mathcal{E},\sum_{[h]\in(\sum^{1}/\sim)\setminus [g]} I_{[h]}\}=0.$$
\noindent It implies $\{x, \mathcal{E}, \mathcal{E} \}=0$. Using the equations $\{I_{[g]}, I_{[h]}, \mathcal{E}\}=\{\mathcal{E}, I_{[g]},  I_{[h]}\}=0$ for  $[g]\neq [h]$, one gets $\{\mathcal{E}, x,  \mathcal{E} \}=\{\mathcal{E},  \mathcal{E},x \}=0$. That is, $x \in \mathrm{Ann}(\mathcal{E})=0$. Thus $x=0$.
\epf

\noindent {\bf Acknowledgements}\quad The authors would like to
thank the referee for valuable comments and suggestions on this
article.

\end{document}